\documentclass{article}
\usepackage{float}
\usepackage{graphicx}
\usepackage{amsmath,amssymb,amsfonts}
\usepackage{mathtools, bm}
\usepackage{amsthm}
\usepackage{url}

\theoremstyle{definition}

\newtheorem{theorem}{Theorem}[section]
\newtheorem{lemma}{Lemma}[section]
\newtheorem{proposition}{Proposition}[section]
\newtheorem{corollary}{Corollary}[section]
\newtheorem{definition}{Definition}[section]
\newtheorem{example}{Example}[section]

\title{A Two-Operator Calculus for Arithmetic-Progression Paths in the Collatz Graph}
\author{\textsc{S. Angermund}}

\begin{document}

\maketitle

\begin{abstract}
\noindent
A recast of the standard residue-class analysis of the \(3x+1\) (Collatz) map in
terms of two elementary operators on arithmetic progressions.  The resulting
calculus \emph{(i)} splits \emph{any} progression into its even and odd
subsequences in a single step; \emph{(ii)} gives a closed formula for every set
of seeds that realises a prescribed parity word; \emph{(iii)} yields a
one–line affine invariant that forbids trajectories consisting of infinitely
many “odd“ moves; and \emph{(iv)} reduces the non-trivial-cycle problem to a pair of linear congruences.
\end{abstract}

\section*{Related Work}\label{sec:related}

The Collatz function has been studied from multiple, partly overlapping,
perspectives.  We outline those most relevant to the present paper.

\paragraph{Stopping-time and density results.}
Terras pioneered the study of stopping times and proved that the set of
integers with finite stopping time has asymptotic density one
\cite{Terras1976}.  Everett subsequently analysed parity sequences to obtain
upper bounds on total stopping times \cite{Everett1977}.  Lagarias gives a
thorough survey and annotated bibliography covering these and many later
developments \cite{LagariasBib2003}.

\paragraph{Residue-class-wise affine (rcwa) mappings.}
A modern unifying language for Collatz-type maps views them as actions of a
subgroup of the rcwa group on \(\mathbb Z\).  Kohl’s monograph
\cite{KohlRCWA2012} formalises this viewpoint and describes algorithmic tools
for rcwa dynamics.  Our operator pair \((T_1,T_2)\) can be seen as generators
of a particularly simple rcwa subgroup acting modulo~\(2\).

\paragraph{Parity words and \(2\)-adic conjugacy.}
The binary “parity word’’ encoding of trajectories dates back to Terras,
reappears in Everett’s parity-sequence bounds, and underlies Akin’s \(2\)-adic
conjugacy map \cite{Akin2004}\cite{LagariasOverview2021}.  We extend
this idea by giving a closed-form expression for the \emph{number} of seeds
realising any finite word, thereby quantifying population density along each
branch of the inverse tree.

\paragraph{Inverse trees and graph-theoretic approaches.}
Numerous authors build a directed tree whose vertices are pre-images of~\(1\);
see, e.g.\ Bernstein–Lagarias’ decomposition of the \(3x+1\) digraph and the
graph-theoretic study by Sultanow \cite{Sultanow2020}.  The tree presented
here differs in that it keeps both affine parameters \((a,b)\) symbolic at
every level, allowing us to track linear invariants such as \(b=a-1\).

\paragraph{Growth-rate versus invariant arguments.}
Classical proofs rule out infinite “all-odd’’ branches by comparing
geometric growth factors (\(3/2\) versus \(1/2\)).  Our proof replaces this
numerical estimate with a single affine constraint \(b=a-1\), providing what
appears to be a new, shorter argument.

\medskip
\noindent
In summary, while our work overlaps with well-established tools—parity sequences,
rcwa dynamics, and inverse trees—it contributes (i) a minimal two-operator
calculus, (ii) an explicit seed-count formula for arbitrary parity words, and
(iii) a new linear-invariant proof excluding trajectories with infinitely many
odd moves.

\section{Terminology}

\begin{definition}
For shorthand we will let $S(a,b)$ refer to the sequence $ (ai+b)_{i\in\mathbb{N}}$ where $a,b\in \mathbb{N}$ and $i$ is the sequence index.
\end{definition}

\section{Divisibility}
Later on we will need to know the divisibility of $s_i\in S(a,b)$ based on $a$ and $b$. With the risk of stating the obvious, let's list all possible cases:
\\[.4cm]
$\forall s_i \in S(a,b)$ it is the case that, when
\begin{enumerate}
    \item
    both $a,b$ even $ \implies 2\mid s_i$,
    \item
    \begin{enumerate}
        \item
        both $a,b$ odd$\implies 2\mid s_i$ for all odd $i$,
        \item
        both $a,b$ odd$\implies 2\nmid s_i$ for all even $i$,
    \end{enumerate}
    \item
    $a$ even and $b$ odd $\implies 2\nmid s_i$,
    \item
    \begin{enumerate}
        \item
        $a$ odd and $b$ even $\implies 2\mid s_i$ for all even $i$,
        \item
        $a$ odd and $b$ even $\implies 2\nmid s_i$ for all odd $i$,
    \end{enumerate}
\end{enumerate}
These are all easy to verify.

\section{Index transformations}
It's interesting for us to consider two specific transformations of $i$, let's define them as $T_1: i \longmapsto 2j$ and $T_2: i \longmapsto 2j+1$. Now $T_1: (ai+b)_{i\in\mathbb{N}} \longmapsto (2aj+b)_{j\in\mathbb{N}}$ and $T_2: (ai+b)_{i\in\mathbb{N}} \longmapsto (2aj+a+b)_{j\in\mathbb{N}}$. Or
\begin{equation*}
\begin{split}
    & T_1: S(a,b) \longmapsto S(2a,b), \\
    & T_2: S(a,b) \longmapsto S(2a,b+a).
\end{split}
\end{equation*}

\begin{example}\label{example:ex1}
Let's look at $S(3,2)$. Apply $T_2: S(3,2)\longmapsto S(6,5)$ and choose for example $s_{10}\in S(6,5) = 65$. Now, since $T_2: i \longmapsto 2j + 1$ and we chose $j=10$ then $i=21$ will give us the same sequence value from the original sequence: $s_{21}\in S(3,2) = 65$.
\end{example}
\noindent
Example~\ref{example:ex1} demonstrates how we can track both the index and the parameters $a,b$ when performing transformations on sequences.

\section{First results}
With the divisibility results from above we can now show why the $T_1, T_2$ transformations are interesting:
\begin{proposition}\label{proposition:prop1}
Parity splittig
\begin{enumerate}
    \item When $a,b$ odd then
    \begin{equation*}
    \begin{split}
        & T_1: S(a,b) \longmapsto (s_i \in S(a,b): 2 \nmid s_i)_{i\in N}, \\
        & T_2: S(a,b) \longmapsto (s_i \in S(a,b): 2 \mid s_i)_{i\in N}.
    \end{split}
    \end{equation*}
    \item When $a$ odd and $b$ even then
    \begin{equation*}
    \begin{split}
        & T_1: S(a,b) \longmapsto (s_i \in S(a,b): 2 \mid s_i)_{i\in N} \\
        & T_2: S(a,b) \longmapsto (s_i \in S(a,b): 2 \nmid s_i)_{i\in N}.
    \end{split}
    \end{equation*}
\end{enumerate}
\end{proposition}
\noindent
For the cases when $a,b$ are both even and when $a$ is even $b$ is odd the sequences are already even or odd.
\\[.4cm]
A proof is not necessary - each case is easily verifiable by looking at the transformed versions $S(2a,b)$ and $S(2a,b+a)$ and using the divisibility results from above. E.g. when $a$ is odd and $b$ is even, then $2a$ is even and $a+b$ is odd. Then according to the divisibility result 3 above - $S(2a,b+a)$ is always odd.
\\[.4cm]
We now have a way to pick out the complete and ordered, even and odd sub-sequences from any sequence $S(a,b)$.

\section{Collatz}
\textsc{The Collatz Conjecture:} For any $n_0\in \mathbb{N}$ the sequence $n_{i+1} = f_{col}(n_i)$ where
\begin{equation*}
    f_{col}(n) = \begin{cases}
        n/2, \quad &when \ 2\mid n\\
        3n+1, \quad &when \ 2\nmid n
    \end{cases}
\end{equation*}
will reach the number 1.
\\[.4cm]
When $n$ is odd, $3n+1 \longrightarrow 3(2\lambda+1)+1 = 6\lambda+4$ is always even. So $f(n)$ can take on a more compact form if we combine the odd-even step. Let's define this function for future reference.
\begin{definition}
The compact Collatz step
\begin{equation}\label{eq:eq1}
    n_{i+1} = \begin{cases}
        n_i/2, \quad &when \ 2\mid n_i\\
        \frac{3n_i+1}{2}, \quad &when \ 2\nmid n_i
    \end{cases}
\end{equation}
\end{definition}

\section{Collatz possibility graph}
Let's start a Collatz sequence without making any assumptions on the seed $n_0$, and build up a graph of possible paths. Making no assumptions on $n_0$ we let $n_0=s_i\in S(1,0)$ be any number (note that $n_0$ is just the sequence index $i$). To perform a Collatz step on $n_0$, we can without loss of generality first do an even and odd sub-sequence extraction: using Proposition~\ref{proposition:prop1}, $S(1,0)$ will branch of into the even sub-sequence $S(2,0)$ and the odd sub-sequence $S(2,1)$.
\begin{enumerate}
    \item
    Even sub-sequence
    \begin{itemize}
        \item
        \textbf{Sub sequence extraction}: $T_1: S(1,0) \longmapsto S(2,0)$
        \item \textbf{Collatz step:} $S(2,0) \longrightarrow S(1,0)$ - \textit{The child loops back to the parent}
    \end{itemize}
    \item
    Odd sub-sequence
    \begin{itemize}
        \item
        \textbf{Sub sequence extraction}: $T_2: S(1,0) \longmapsto S(2,1)$
        \item \textbf{Collatz step:} $S(2,1) \longrightarrow S(6,4) \longrightarrow S(3,2)$
    \end{itemize}
\end{enumerate}
In the Collatz steps we have simply plugged in $n_i=ai+b$ into Equation~\ref{eq:eq1}.
\\[.4cm]
\textsc{Remark 1:} Without loss of generality we have here constructed the \textit{form} that all possible numbers $n_1 = f_{col}(n_0)$ can take on. Using Proposition~\ref{proposition:prop1} we can keep applying this method to the child nodes, and without loss of generality we can construct the \textit{structure} $S(a,b)$ of all nodes in all possible paths of Equation~\ref{eq:eq1}.
\\[.4cm]
\textsc{Remark 2:} While sub-sequence extractions will transform $i$, the Collatz steps will preserve $i$. Even: $f_{col}(ai+b) = (a/2)i+(b/2)$. Odd: $f_{col}(ai+b)=3ai + (3b+1)$.

\subsection*{Tracing the seed}
A sequence of transformations $T_{k\in K}$, where $K$ is a sequence of 1´s and 2´s, has the following property:
\begin{equation*}
    T_{k\in K}: i \longmapsto 2^{|K|}j+x
\end{equation*}
where $x\geq0$ is some sum of powers of 2. More generally:
\begin{proposition}[General index formula]\label{proposition:prop2}
Let $K=(\kappa_{1},\dots ,\kappa_{k})\in\{T_{1},T_{2}\}^{k}$ be a word in the
two operators, read from root to leaf, and write
\[
s_{t}\;=\;
  \begin{cases}
    1 &\text{if } \kappa_{t}=T_{2},\\
    0 &\text{if } \kappa_{t}=T_{1}.
  \end{cases}
\]
After the $k$ subsequence extractions the seed index $i$ that lands on
position $j$ in the target progression is
\[
\boxed{\;
  i \;=\; 2^{\,k}\,j \;+\; \sum_{t=1}^{k} s_{t}\,2^{\,k-t}\; }.
\]
\noindent
\textsc{Proof:} Each operator doubles the running index; whenever $\kappa_{t}=T_{2}$ it
adds~$1$.  Unrolling the recursion gives
\(
i = 2^{k} j
    + s_{1} 2^{k-1}
    + s_{2} 2^{k-2}
    + \cdots
    + s_{k}.
\)
\end{proposition}
\begin{flushright}
$\Box$
\end{flushright}
So when only $T_1$ transformations are applied then $T_{k\in K}:i\longmapsto 2^{|K|}j, \ K=(1,1,1,...,1)$. When only $T_2$ transformations are applied we will get the final transformation $T_{k\in K}:i\longmapsto 2^{|K|}j+2^{|K|}-1, \ K=(2,2,2,...,2)$. This follows from the fact that $2^{k-1}+...+2^0=2^k-1$.
\\[.4cm]
We can now trace the seed in the graph by tracking the transformation sequence $K$ on a given path. There is also an interesting result following this regarding path densities:

\begin{corollary}[Seed density of a node]
For a word $K$ of length $k$, the set of seeds that reach the
corresponding node is the arithmetic progression
$2^{k}\mathbb{N}+\sum_{t=1}^{k}s_{t}2^{k-t}$.
Consequently the sub-tree rooted at that node has natural
density $2^{-k}$ inside $\mathbb N$.
\end{corollary}

\subsection*{Illustrations}
Running the first few steps on a computer we can produce some visualizations. Here is the first few levels of the tree structure that forms when running the compact Collatz function:
\begin{figure}[H]
\begin{center} 
\includegraphics[width=1.0\textwidth]{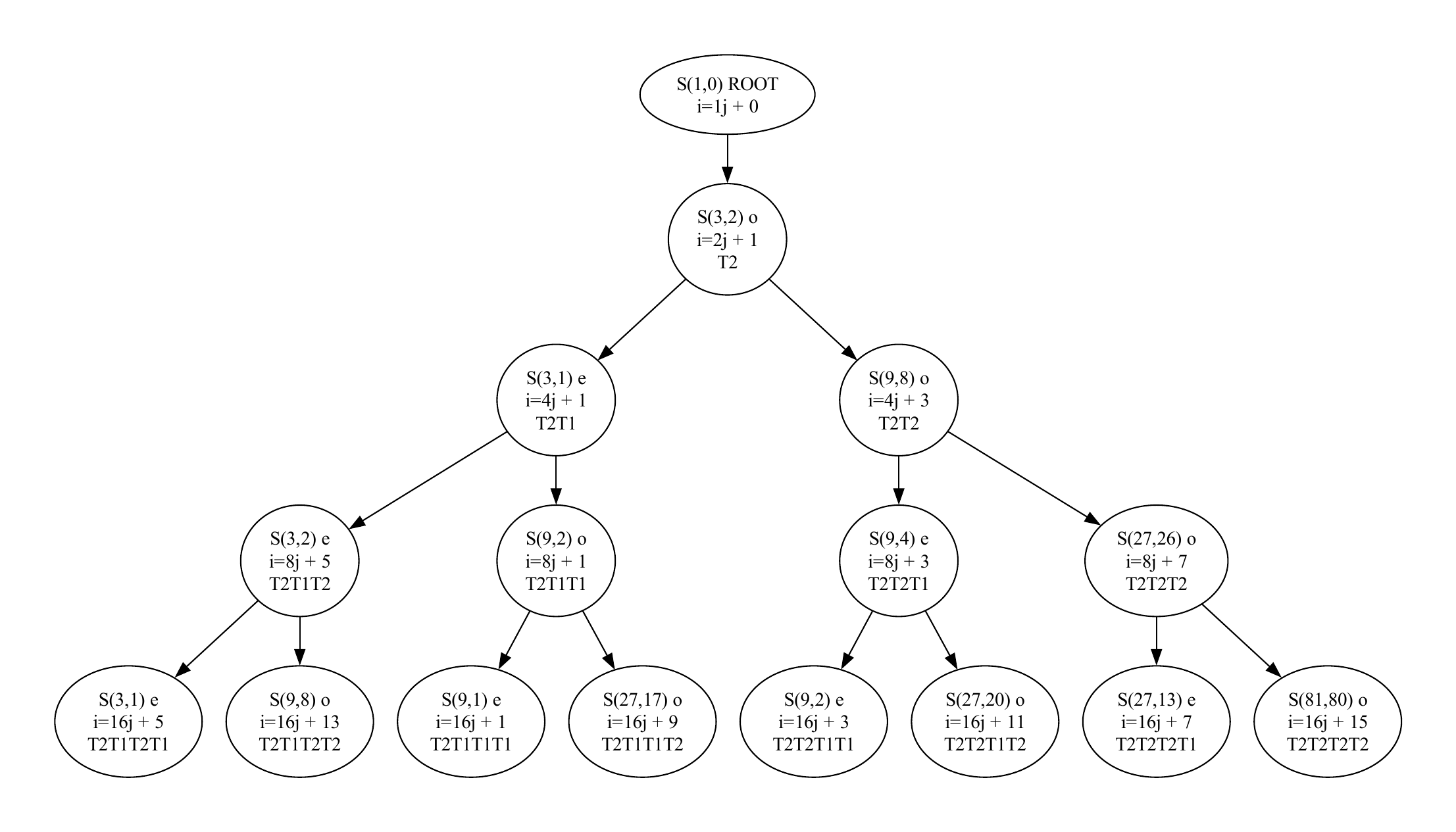}
\caption{A visualization of the compact Collatz steps for any seed. Each node shows the arithmetic sequence S(a,b); prefix ‘even/odd parent’ records whether the incoming Collatz step used the even or odd branch. ‘i = …’ shows the accumulated index map and below the word in $\{T_1,T_2\}$ to reach the node.}
\label{fig:fig2}
\end{center}
\end{figure}
Link to the python code used to generate the tree object and visualizations is referenced \cite{Angermund2025}
\begin{example}\label{example:ex2}
Let's look at Figure~\ref{fig:fig2} and trace the lineage of $i=n_0$ down to the lower left most node $S(9,5)$. Remember that we have only done a sub-sequence extraction when both $a,b$ are odd or $a$ odd and $b$ even in accordance to Proposition~\ref{proposition:prop1}:
The first transformation is the odd sub-sequence step from $S(1,0)$ which transforms the index as $T_2$, then the even sub-sequence from $S(3,2)$ which transforms the index as $T_1$, and so on. Resulting in the transformation sequence $T_2T_1T_1T_1T_2(i)$. I.e. $i \longrightarrow 2(2\cdot2\cdot2\cdot(2j+1))+1 = 2(16j+8)+1 = 32j+17$, where $j$ is the index of the target node $S(9,5)$.
So, by choosing a $j$, we will select both the node value and the seed value. Let's pick $j=5 \implies s_5\in S(9,5) = 50$ and $n_0 = i = 32j+17 = 177$. We verify this by running the Collatz algorithm with seed 177 and checking that the 5th compact (7th regular) step is 50. The resulting values are 177, 266, 133, 200, 100, 50.
\end{example}
Example~\ref{example:ex2} showcases that, since the value of $j$ was arbitrary, there are infinitely many seeds that can lead to a given sequence of Collatz steps (a path in the graph). In the next section we will look at what this can teach us about general Collatz sequences.
\section{Odd steps only}
A not so far fetched guess would be that there is a seed $n_0$ such that the odd step of the compact Collatz function would be infinitely repeated: successions of multiplying by 3, adding 1 and dividing by 2 will clearly diverge to infinity as long as $n_0\geq0$. We will now attempt to show that there is no such seed.
\begin{theorem}\label{theorem:thm1}
There is no finite seed $n_0$ such that Equation~\ref{eq:eq1} will repeat the odd step indefinitely.
\\[.2cm]
\textsc{Proof:} We can assume that $n_0$ is odd and thereafter the Collatz step will always be
\begin{equation*}
    n_{k+1} = \frac{3n_k + 1}{2}.
\end{equation*}
This is the right most path in Figure~\ref{fig:fig2}, and by a glance it looks like $n_k$ always has an odd $a$ and even $b$, in fact the pattern seems to be that $a$ is odd and $b=a-1$. If we can show that this is always the case then, according to Proposition~\ref{proposition:prop1}, we are always doing $T_2$ index transformations. Let's prove this as a lemma:
\begin{lemma}
We want to show that for $n_k = ai+b$, where $a$ odd and $b=a-1$, after a $T_2$ transformation and the odd step in Equation~\ref{eq:eq1} we will end up with $n_{k+1}=a'j+b'$ where $a'$ is odd and $b'=a'-1$. Let's do the full induction proof.
\\[.2cm]
\textbf{Base case:} $n_0 = s_i \in S(1,0)$. Odd sub-sequence extraction: $T_2: S(1,0) \longrightarrow S(2,1)$. Next we do the Collatz step on $s_i\in S(2,1)$:
\begin{equation*}
    n_1 = \frac{3(2i+1)+1}{2} = \frac{6i+4}{2} = 3i+2
\end{equation*}
3 is indeed odd and $2=3-1$ so base case complete.
\\[.2cm]
\textbf{Induction case:} $n_k = s_i \in S(a,b)$. Assume $a$ odd and $b=a-1$. Odd sub-sequence extraction: $T_2: S(a,b) \longmapsto S(2a, b+a)$, and Collatz step on $s_i \in S(2a,b+a)$:
\begin{equation*}
    n_{k+1} = \frac{3(2ai+b+a)+1}{2} = \frac{6ai+3b+3a+1}{2} = 3ai + \frac{3(a+b)+1}{2}
\end{equation*}
We can apply the transformations $a=2\lambda+1$ and $b=2\lambda$ as the induction hypothesis:
\begin{equation*}
\begin{split}
    n_{k+1} &= 3ai + \frac{3(a+b)+1}{2} = 3(2\lambda+1)i + \frac{3(2\lambda+1+2\lambda)+1}{2}\\
    &= (6\lambda+3)i + \frac{12\lambda+4}{2} = (6\lambda+3)i + (6\lambda+2)
\end{split}
\end{equation*}
Now $a'=6\lambda+3$ and $b'=6\lambda+2$, which gives the resulting sequence $S(a',b')$ where indeed $a'$ is odd and $b'=a'-1$.
\end{lemma}
\begin{flushright}
$\Box$
\end{flushright}
Now, importantly, we know that we are always doing $T_2$ index transformations before each Collatz step. So after $k$ successive steps we will land at a node with sequence index $j$, and by~\ref{proposition:prop2} we can retrieve the seed as $s_i\in S(1,0) = i = 2^kj + 2^k -1$. Even if we choose $j=0$, the starting seed will clearly be unbounded when we let $k\longrightarrow\infty$. Thus there is no finite seed $n_0$ such that the compact odd Collatz-step can be applied indefinitely.
\begin{flushright}
$\Box$
\end{flushright}
\end{theorem}

\section{Cycles}\label{sec:cycles}
We can use the parity splitting together with the index tracing to search for cycles.

\begin{proposition}[Cycle conditions]\label{proposition:prop3}
Let $n_k$ be a node at depth $k$ with arithmetic sequence $a_1 j + b_1$ and index map $\alpha_1 j + \beta_1$. Let $n_{k-t}$ be a node at depth $k-t$ with arithmetic sequence $a_2 j + b_2$ and index map $\alpha_2 j + \beta_2$.  For there to be a cycle between $n_k$ and its ancestor $n_{k-t}$ there must exist $x,y\in\mathbb{N}$ such that
\begin{equation*}
    \begin{split}
        a_1 x + b_1 & = a_2 y + b_2, \\
        \alpha_1 x + \beta_1 & = \alpha_2 y + \beta_2.
    \end{split}
\end{equation*}
This ensures that the nodes take on equal value when spawned from the same seed.
\end{proposition}
If we solve the linear system for $x$ and $y$ we get an explicit formula for cycle conditions:
\begin{equation*}
    \begin{split}
        x & = \frac{a_2\big(\beta_1 - \beta_2\big) + \alpha_2\big(b_2 - b_1\big)}{a_1\alpha_2 - a_2\alpha_1}, \\
        y & = \frac{a_1\big(\beta_1 - \beta_2\big) + \alpha_1\big(b_2 - b_1\big)}{a_1\alpha_2 - a_2\alpha_1}.
    \end{split}
\end{equation*}
Notice that if $b_1=b_2$ and $\beta_1=\beta_2$ we get $x=y=0$. We can actually see this cycle in the first few levels of the tree graph.
\begin{figure}[H]
\begin{center} 
\includegraphics[width=1.0\textwidth]{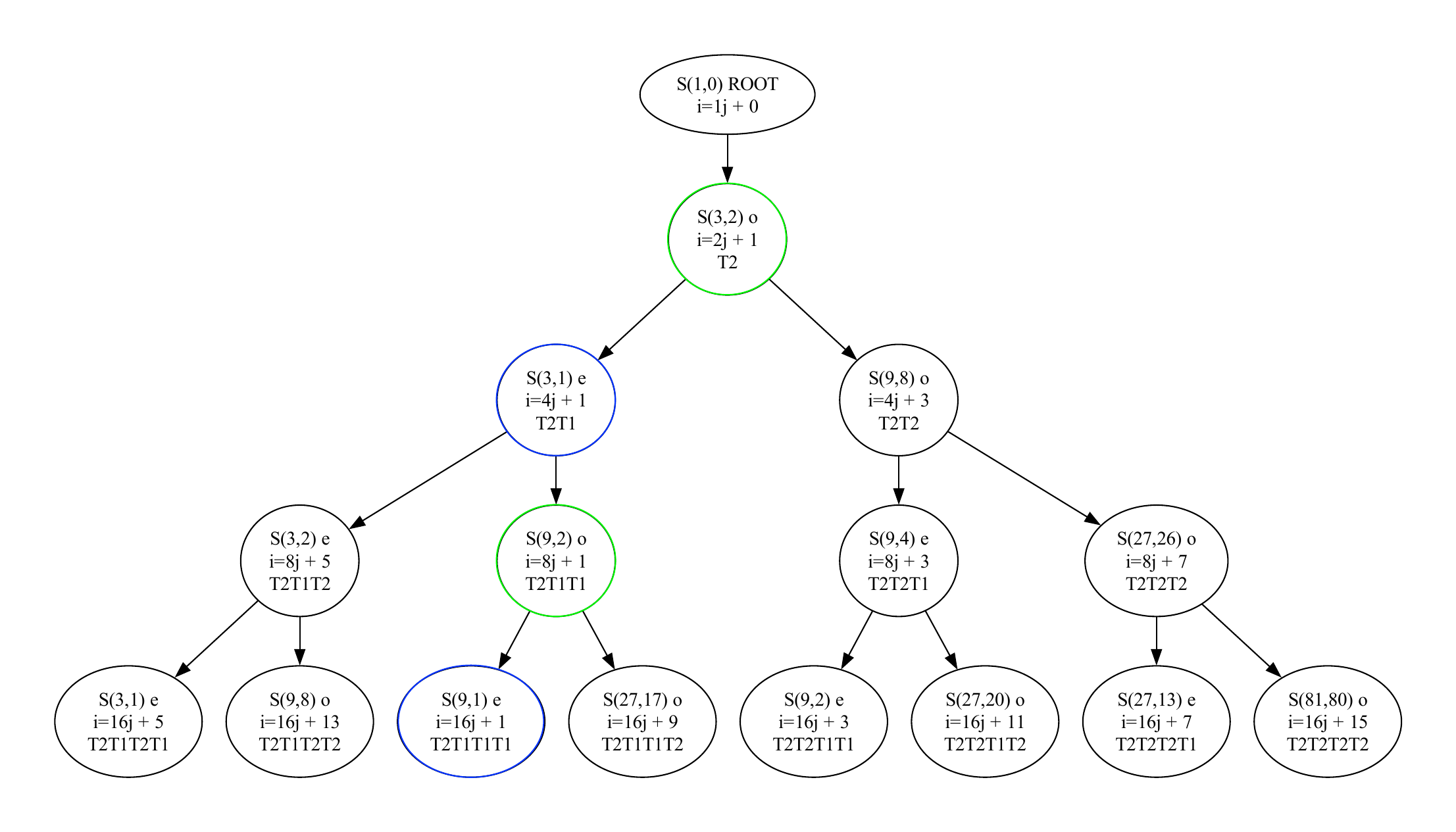}
\caption{If we let $j=0$ on the $S(3,2)$ node in level 1 you get node value 2 and seed 1. The same for the $S(9,2)$ node in level 3. These are examples of the 1-2-1-2.. cycle of the compact Collatz function. Note that these nodes satisfy the $b_1=b_2,\beta_1=\beta_2$ condition.}
\label{fig:fig3}
\end{center}
\end{figure}
Rudimentary code to find cycles is referenced \cite{Angermund2025}. Whether this framework can be useful in contrast to methods such as the m-cycle matrix technique\cite{SimonsDeWeger2005} is unclear to us at this stage.

\section{Conclusion and Outlook}\label{sec:conclusion}

We have shown that the Collatz dynamics admit a compact
description once every integer is encoded as an arithmetic progression
\(S(a,b)\).  A single pair of index operators,
\[
  T_{1}:S(a,b)\longmapsto S(2a,b),
  \qquad
  T_{2}:S(a,b)\longmapsto S(2a,b+a),
\]
\emph{(i)} separates any progression into its even and odd
subsequences in one step; \emph{(ii)} yields a closed formula
\(i\mapsto 2^{k}j+2^{k}-1\) for the seeds that realise any prescribed
parity word (Proposition~\ref{proposition:prop2}); and \emph{(iii)} exposes the
affine invariant \(b=a-1\) that rules out trajectories consisting of
infinitely many odd moves (Theorem~\ref{theorem:thm1}).

\subsection{Immediate corollaries.}
\begin{itemize}
  \item Every Collatz orbit must include at least one division by \(2\); the
        invariant argument replaces the classical growth-rate estimate
        \(3/2\) versus \(1/2\).
  \item Each node in the symbolic inverse tree carries an explicit
        population count—the number of seeds that flow through
        it—clarifying why some branches are denser than others.
\end{itemize}

\medskip
\noindent
\subsection{Directions for further work.}
\begin{enumerate}
  \item \emph{Cycle classification.}  
        We started analysing cycles in section \ref{sec:cycles}; a systematic search within the
        \(\{T_{1},T_{2}\}\) calculus could sharpen existing cycle-length
        bounds.
  \item \emph{Density estimates.}  
        The seed-count formula suggests a probabilistic model for the
        likelihood that a random integer follows a given parity word.
  \item \emph{Generalised maps.}  
        Replacing the pair \((3n+1,n/2)\) by other residue-class-wise
        affine rules (e.g.\ \(5n+1\) on odd numbers) preserves the
        operator framework; studying which invariants survive such
        changes may shed light on why \(3x+1\) is special.
  \item \emph{Algorithmic applications.}  
        Because the operators act symbolically on \((a,b)\), they can
        serve as a state-space reduction in exhaustive searches or SAT
        encodings of the conjecture.
\end{enumerate}

\bibliographystyle{amsplain}   % or plain, alpha, abbrv, etc.
\bibliography{references}      % name without the .bib extension

\end{document}